\documentclass[12pt]{article}
\usepackage{amsfonts}
\usepackage{wasysym,amssymb,eufrak,indentfirst,graphicx,cite,amsthm,color}
\usepackage[bookmarksnumbered, colorlinks, plainpages]{hyperref}
\newtheorem{Theorem}{Theorem}[section]
\newtheorem{Corollary}[Theorem]{Corollary}
\newtheorem{Lemma}[Theorem]{Lemma}
\newtheorem{Proposition}[Theorem]{Proposition}
\newtheorem{Question}[Theorem]{Question}
\newtheorem{Definition}[Theorem]{Definition}

\newtheorem{Remark}{{Remark}}

\def\qed{\hfill $\Box$}

\begin{document}
\baselineskip 17pt

\title{On nilponent hypergroups\thanks{Research was supported by the NNSF  of China (12001526, 12171126)
and Natural Science Foundation of Jiangsu Province, China (BK20200626).}}

\author{Chi Zhang\\
{\small  Department of Mathematics, China University of Mining and
Technology}\\
{\small Xuzhou, 221116, P. R. China}\\
{\small E-mail: zclqq32@cumt.edu.cn}\\ \\
Wenbin Guo\thanks{Corresponding author}\\
{\small School of Science, Hainan University}\\
{\small  Haikou, Hainan, 570228, P.R. China, \ and }\\
{\small School of Mathematical Sciences, University of Science and Technology of China}\\
{\small Hefei 230026, P.R. China}\\
{\small E-mail: wbguo@ustc.edu.cn}\\}

\date{}
\maketitle

\begin{abstract}
In this paper, we establish the theory of nilpotent hypergroups
and study some properties of nilpotent hypergroups and provided some structural characterizations of nilpotent hypergroups.
\end{abstract}

\let\thefootnoteorig\thefootnote
\renewcommand{\thefootnote}{\empty}


 \let\thefootnote\thefootnoteorig

\section{Introduction}
F. Marty \cite{m1} introduced the notion of hypergroups which generalizes the definition of groups.
The following definitions regarding hypergroups are due to P.H.-Zieschang and various co-authors \cite{b1,fz1,vz1}.

\begin{Definition}

A hypergroup is a set $H$ equipped with a hypermultiplication (a map from $H$ to the power set of $H$, denoted as $(p, q)  \longmapsto pq$ for all $p, q \in H$).
For $P, Q$ any subsets of $H$,
$$PQ := \bigcup_{p \in P, q \in Q} pq.$$
If $P = \{ p \}$, a singleton set, then $pQ := \{ p \}Q$, and $Qp := Q \{ p \}$.
The hypermultiplication is assumed to satisfy the following conditions:

$(H1)$  For any elements $p, q$, and $r$ in $H$, $p(qr) = (pq)r$.

$(H2)$ $H$ contains an element $1$ such that $s1 = \{ s \}$ for all $s \in H$.

$(H3)$ For each element $s$ in $H$, there exists an element $s^*$ in $H$ such that for any elements $p, q$, and $r$ in $H$ with $r \in pq$, then we have $q \in p^*r$ and $p \in rq^*$.

\end{Definition}

$(H1)$ implies that set product is associative. $(H2)$ and $(H3)$ yield that for all $s \in H$, $1 \in s^*s$.
An element $s \in H$ is called {\em thin} if $s^*s = \{ 1 \}$.

The set of all the thin elements of a hypergroup $H$ is denoted $O_{\vartheta}(H)$.
A hypergroup $H$ is called thin if $H = O_{\vartheta}(H)$, that is, all elements of $G$ are thin.
Of course, any group $H$ may be regarded as a thin hypergroup, simply by replacing the product of two elements with the singleton set containing that product.
Conversely, it is easy to see that a thin hypergroup is a group.
In fact, the concept of hypergroups is a generalization of the concept of groups.

Later, a number of group theoretic results have found a generalization within the theory of hypergroups.
In \cite{tz1}, the Homomorphism Theorem and the First Isomorphism Theorem for groups were generalized to structure theorems on hypergroups by R. Tanaka and P.-H. Zieschang.
In \cite{z2}, P.-H. Zieschang has shown that quite a few facts on involutions can be generalized from group theory to the theory of hypergroups.
In \cite{b1}, H. Blau generalizes the Sylow Theorem of finite group to the hypergroup.
In 2022, A. V. Vasil'ev and P.-H. Zieschang \cite{vz1} proposed the definition of solvable hypergroups and generalized the theory of finite solvable groups to solvable hypergroups.
As we known,  solvable groups and nilpotent groups are both impotent classes on the theory of groups.
So we naturally propose the following question.

\vspace{0.15cm}

{\bf Question.} {\sl Can we establish the theory of nilpotent hypergroups ?}

\vspace{0.15cm}

In this paper, we will answer this above question(see below Definition \ref{nilpotent}).

Following \cite{fz1, vz1},
for any subset $S$ of hypergroup $H$, let $S^* := \{ s^* | s \in S \}$;
A nonempty subset $F$ of $H$ is called {\em closed} if $F^*F \subseteq F$, that is, $a^*b \subseteq F$ for all $a, b \in F$.
Clearly, the intersection of some closed subsets of $H$ is also a closed subset of $H$.
A closed subset $F$ of $H$ is called {\em normal} in $H$ if $Fh \subseteq hF$ for each element $h$ in $H$;
A closed subset $F$ of $H$ is called {\em strongly normal} in $H$ if, for each element $h$ in $H$, $h^*Fh \subseteq F$.
It is easy to see that strong normality implies normality.

The intersection of all closed subsets of $H$ which are strongly normal in $H$ will be denoted by $O^{\vartheta}(H)$.
The subset $O^{\vartheta}(H)$ of $H$ is called the {\em thin residue} of $H$.
Lemma $3.6(ii)$ in \cite{fz1} implies that the thin residue of $H$ is a strongly normal closed subset of $H$.

Following \cite{fz1, vz1},
a closed subset $F$ is {\em subnormal} ({\em strongly subnormal}) in $H$ if
there exists a chain of closed subsets $F = F_0 \subseteq F_1 \subseteq \cdots \subseteq F_n = H$ for some $n > 0$ such that $F_{i-1}$ is normal (strongly normal) in $F_i$ for all $1 \leq i \leq n$.
Let $F$ be a closed subset of a hypergroup $H$.
For any $h \in H$, define $h^F := FhF$; and for any subset $F \subseteq S \subseteq H$, define $S//F := \{ h^F | h \in S \}$.
Then for all $a, b \in H$, $a^F \cdot b^F := \{x^F | x \in aFb \}$
defines a hypermultiplication on $H//F$ such that $H//F$ becomes a hypergroup, called the quotient of $H$ over $F$.
Recall that $[a, b]=a^*b^*ab$ and for any two subsets $A$ and $B$ of $H$,  $[A, B]$ denotes the smallest closed subset of $H$ which contains all sets $[a, b]$ with $a \in A$ and $b \in B$.
It is easy to see that $[A, B]$ is the intersection of all closed subsets of $H$ containing all sets $[a, b]$  with $a \in A$ and $b \in B$.

\begin{Definition}\label{nilpotent}

A hypergroup $H$ is said to be nilpotent if $H_{n} = \underbrace{[H, H, \cdots, H]}_{n} = 1$ for some positive integer $n$.
Where $H_{n} = [H_{n-1}, H]$ and $H_{1} = H$.

\end{Definition}

We will study the characterisetion of nilpotent hypergroups and get some results that generalize the classic results of nilpotent groups.
Let $Z(H) = \{ h \in H | hx=xh, \forall x \in H \}$ and we call that $Z(H)$ is the centre of $H$.
However in general, $Z(H)$ is not closed in $H$.
We define that $Z^*(H) = \{ x \in Z(H) | x^* \in Z(H) \}$.
It is easy to see that $Z^*(H)$ is a normal closed subset of $H$ (see below Lemma \ref{main1}).
We call that $Z^*(H)$ is a closed center of $H$.
Then we give the following definition.

\begin{Definition}\label{hyper}

Let $H$ be a hypergroup.
We call that a closed subset series of $H$: $$1 = Z_{0}^{*}(H) \subseteq Z_{1}^*(H) \subseteq \cdots \subseteq Z_{n}^*(H) \subseteq \cdots$$  the closed center series of $H$
if $Z_{i}^*(H)//Z_{i-1}^*(H) = Z^*(H//Z_{i-1}^*(H))$ for all $i = 1, 2, \cdots$.

\end{Definition}

\begin{Remark}

We could show that $Z_{i}^*(H)$ is a normal closed subset of $H$ for all $i = 1, 2, \cdots$ (See Lemma \ref{main1}).
When $Z_{n}^*(H) = Z_{n+1}^*(H) = \cdots$, then we denote by $Z_{\infty}^*(H)$ the terminal term of this ascending series and say that $Z_{\infty}^*(H)$ is the {\em inv-hypercenter} of $H$.

\end{Remark}

We get the following theorem.

\begin{Theorem}\label{center}

If A hypergroup $H$ is nilpotent, then $Z_{n}^*(H) = H$ for some non-negative integer $n$.

\end{Theorem}

And we get a sufficient condition for judging the hypergroup $H//O^{\vartheta}(H)$  to be nilpotent.

\begin{Theorem}\label{ct}

If $Z_{n}^{*}(H) = H$ for some non-negative integer $n$, then $H//O^{\vartheta}(H)$ is nilpotent.

\end{Theorem}

We also get some characterization of finite nilpotent hypergroups.

\begin{Theorem}\label{strongly}

Let $H$ be a finite hypergroup. If $H$ is nilpotent, then every closed subset $1 \neq E$ of $H$ is strongly subnormal in $H$.

\end{Theorem}

\begin{Remark}
From Theorem \ref{strongly}, we get that if $H$ be a finite nilpotent hypergroup, then every closed subset $E \neq 1$ of $H$ is strongly subnormal in $H$.
In fact, $1\cdot h=h\cdot 1$ for all $h \in H$, so $1$ is normal in $H$.
However, $h^{*}1h \subseteq 1$ is not true in general.
\end{Remark}

As we known, every nilpotent group is a soluble group in the theory of groups.
A hypergroup $H$ is said to be {\em solvable} \cite{vz1,b1} if it contains closed subsets $F_{0}, \cdots, F_n$ such
that $F_0 = {1}$, $F_n = H$, and, for each element $i$ in $\{0, \cdots, n\}$ with $1 \leq i$, $F_{i-1} \leq F_i$,
$F_i//F_{i-1}$ is thin, and $|F_i//F_{i-1}|$ is a prime number.
Now in the theory of finite hypergroups, we get the following theorem.

\begin{Theorem}\label{ns}

If $H$ is a finite nilpotent hypergroup, then $H$ is a solvable hypergroup.

\end{Theorem}

This paper is organized as follows. In section $2$, we cite some known results which are useful
in our proofs and prove some basic properties of hypergroups.
In section $3$, we study the nilpotent hyprergroups, give the proofs of Theorems \ref{center}, \ref{ct}, \ref{strongly} and \ref{ns}
and also give some charaterizations of nilpotent hypergroups.
In section $4$, we give some remarks and questions on nilpotent hypergroup theory.

In this paper, the letter $H$ stands for a  hypergroup.

\section{Some known results and basic preliminaries}

\begin{Lemma}\label{basic}{\rm (See \cite[Lemma 2.1, Lemma 3.1(ii) and Corollarly 2.6]{fz1})}
Let $a, b$ be  elements in $H$.
Let $A, B$ be subsets of $H$.
Then the following hold.

$(1)$  $a^{**} = a$, where $a^{**} = {(a^{*})}^{*}$

$(2)$ $(ab)^* = b^*a^*$.

$(3)$ $(AB)^* = B^*A^*$.

\end{Lemma}

\begin{Lemma}\label{inv}

Let $p$ and $q$ be elements in $H$. Then the following hold:

$(1)$ $1 \in pq$ if and only if $p=q^*$($q=p^*$) (See \cite[Lemma 1.1(iii)]{tz1}).

$(2)$ $1s=s$. $H$ possesses exactly one neutral element $1$ (See \cite[Lemma 1.3(1)]{tz1}).

\end{Lemma}

We define $C_H(F) = \{ h \in H | hf=fh, \forall f \in F \}$.
$C_H(F)$ will be call the centralizer of $F$ in $H$.
Note that $E \subseteq C_H(F)$  if and only if $F \subseteq C_H(E)$, and $Z(H) = C_H(H)$.

\begin{Corollary}\label{abel}
Let $a, b$ be  elements in $H$. Then the following assertions are equivalent.

$(1)$ $ab=ba$

$(2)$ $1 \in a^*b^*ab$

$(3)$ $a \in C_H(b)$ and $b \in C_H(a)$

\end{Corollary}

{\bf Proof.}

It is easy to see this corollary follows from Lemma \ref{basic} and \ref{inv}.

\qed

\begin{Lemma}\label{closed}

$(1)$  A subset A of $H$ is closed if and only if $1 \in A$, $A^{*} = A$ and $AA = A$. (see \cite[Lemma 3.3]{fz1})

$(2)$ $Fh = hF$ for each normal closed subset of $H$ and any element $h$ in $H$. (see \cite[Lemma 3.1(iii) and p.99]{fz1})

\end{Lemma}

\begin{Lemma}\label{cen}
Let $H$ and $F$ be a hypergroup.
If $[H, F] = 1$,
then for every $h$ in $H$, $h$ centralizes $F$, and for every element $f$ in $F$, $f$ centralizes $H$.
\end{Lemma}

{\bf Proof.}
If $[H, F] = 1$, then for every $h$ in $H$ and $f$ in $F$, we have $[h, f] =1,$ that is, $h^*f^*hf=1.$
By Corollary \ref{abel}, $hf=fh,$ that is, $h$ centralizes $F$.
Similarly, for every element $f$ in $F$, $f$ centralizes $H$.
\qed
\medskip

\begin{Lemma}\label{thin residue}

Let $H$ be a hypergroup.

(1) $O^{\vartheta}(H) = [H, {1}]$ (See \cite[Lemma 5.4]{fz1}).

(2) $O^{\vartheta}(H)$ is strongly normal in $H$ (See \cite[Lemma 3.6(ii)]{fz1}).

\end{Lemma}

\begin{Lemma}\label{normal}

If $N$ is normal in $H$, then $F//N$ is normal in $H//N$ if and only if $F$ is normal in $H$.

\end{Lemma}

{\bf Proof.}
$F//N$ is normal in $H//N$ if and only if $(NFN)(NhN) \subseteq (NhN)(NFN)$ if and only if $Fh \subseteq hF$ for every $h \in H$ if and only if $F$ is normal in $H$.
\qed

\begin{Lemma}\label{qinv}

Let $F$ be a closed subset of $H$ and $h$ be a element in $H$.
Then ${(h^{F})}^{*} = {(h^{*})}^{F}.$

\end{Lemma}

{\bf Proof.}
By Lemma \ref{basic},
$${(h^{F})}^{*} = {(FhF)}^{*} = F^{*}h^{*}F^{*} = Fh^{*}F.$$
Hence ${(h^{F})}^{*} = {(h^{*})}^{F}.$

\qed

Following \cite[p.100]{fz1}, a subset $A$ of $H$ is said to be a {\em star invariant} subset of $H$ if $A^{*} = A.$
For each subset $A$ of $H$,  $\langle A \rangle$ denote to be the intersection of all closed subsets of
$H$ containing $A$.

\begin{Lemma}\label{de} {\rm (See \cite[Lemma 3.9]{fz1})}

Let $A$ be a star invariant subset of $H$. Then the following hold.

$(i)$ The set $\langle A \rangle$ is equal to the union of the sets $A^{n}$ with $n$ a non-negative integer.

$(ii)$ Assume that $A$ is not empty, and let $h$ be an element of $H$ satisfying $h^{*}Ah \subseteq \langle A \rangle $.
Then $h^{*}\langle A \rangle h \subseteq \langle A \rangle $.

\end{Lemma}

\begin{Lemma}\label{strong} {\rm (See \cite[Lemma 3.7]{vz1})}

A closed subset $F$ of $H$ is strongly normal in $H$ if and only if $H//F$ is thin.

\end{Lemma}

\begin{Lemma}\label{qu} {\rm (See \cite[Lemma 5.5]{fz1})}

Let $F$ be a normal closed subset of $H$. Then $O^{\vartheta}(H)F//F=O^{\vartheta}(H//F)$

\end{Lemma}

\section{Proofs of Theorems  \ref{center}, \ref{ct}, \ref{strongly} and \ref{ns}}

In this section, we will prove Theorems \ref{center}, \ref{ct}, \ref{strongly} and \ref{ns}.
And we will also give some characterizations of nilpotent hypergroups.
To prove our theorems, we first proved the following main lemmas, which are the key steps in proving the theorems.

\begin{Lemma}\label{main1}

$Z_{n}^{*}(H)$ is normal closed subsets of $H$ for $n = 1, 2, \cdots$.

\end{Lemma}

{\bf Proof.}
For $n =1$, $Z_{1}^*(H) = Z^*(H)$ is a normal closed subset.
In fact, for any two $a, b \in Z^*(H)$, we have $a^* \in Z(H)$ and $b \in Z(H)$.
For every element $h$ in $H$, we have $a^*bh = a^*hb =ha^*b$ and so $a^*b \in Z(H)$.
And $(a^*b)^* = b^*a$ by Lemma \ref{basic}(ii).
Since $b \in Z^*(H)$, $b^* \in Z(H)$.
So similarly, we have $(a^*b)^* = b^*a \in Z(H)$.
Therefore $a^*b \in Z^*(H)$, that is, $Z^*(H)$ is a closed subset.
Since $Z^*(H) \subseteq Z(H)$, $Z^*(H)h = hZ^*(H)$.
Hence $Z^*(H)$ is a normal closed subset.
We assume that $Z_{i}^*(H)$ is a normal colsed subset of $H$.
We show that $Z_{i+1}^*(H)$ is a normal closed subset of $H$.
Firstly, we show $Z_{i+1}^*(H)$ is a closed subset of $H$.
For every $x, y \in Z_{i+1}^*(H)$, $x^{Z_{i}^*(H)}, y^{Z_{i}^*(H)} \in Z^*(H//Z_{i}^*(H))$.
Then $(xy)^{Z_{i}^*(H)} \subseteq x^{Z_{i}^*(H)}y^{Z_{i}^*(H)} \in Z^*(H//Z_{i}^*(H))$.
Therefore $xy \in Z_{i+1}^*(H)$.
Also $(x^*)^{Z_{i}^*(H)} = (x^{Z_{i}^*(H)})^{*} \in Z^*(H//Z_{i}^*(H))$ by Lemma \ref{qinv}.
Hence $Z_{i+1}^*(H)$ is a closed subset of $H$.
Finally, we show that $Z_{i+1}^*(H)$ is  normal in $H$.
In fact, as $Z_{i+1}^*(H)//Z_{i}^*(H) = Z^*(H//Z_{i}^*(H))$,
from the above proof we see that $Z_{i+1}^*(H)//Z_{i}^*(H)$ is normal in $H//Z_{i}^*(H)$.
Now by Lemma \ref{normal}, $Z_{i+1}^*(H)$ is normal in $H$.
Thus, we get that $Z_{n}^*(H)$ is a normal closed subset of $H$ for all $n = 1, 2, \cdots$.
\qed
\medskip

\begin{Lemma}\label{cq}

Let $F$ be a normal closed subset of $H$.
And let $C$ and $D$ be a closed subset of $H$.
Then $[C//F, D//F] = [C, D]F//F.$

\end{Lemma}

{\bf Proof.}
Let $x^{F} \in [C//F, D//F]$. Then it implies from Lemma \ref{de}(i) that
$$x^{F} \in [c_1^{F}, d_1^{F}][c_2^{F}, d_2^{F}] \cdots [c_n^{F}, d_n^{F}]=((c_1^{F})^{*}(d_1^{F})^{*}c_1^{F}d_1^{F}) \cdots ((c_n^{F})^{*}(d_n^{F})^{*}c_n^{F}d_n^{F})$$
where $c_i \in C$ and $d_i \in D$ for some non-negative $n$ and $i = 1, 2, \cdots, n.$
Since $F$ is normal in $H$ and by Lemma \ref{qinv},
$$((c_1^{F})^{*}(d_1^{F})^{*}c_1^{F}d_1^{F}) \cdots ((c_n^{F})^{*}(d_n^{F})^{*}c_n^{F}d_n^{F}) = (c_1^{*}d_1^{*}c_1d_1)^{F} \cdots (c_n^{*}d_n^{*}c_nd_n)^{F}.$$
Therefore $x^{F} \in (c_1^{*}d_1^{*}c_1d_1)^{F} \cdots (c_n^{*}d_n^{*}c_nd_n)^{F}$
and so $x^{F} \in [C, D]F//F.$
It follows that $[C//F, D//F] \subseteq [C, D]F//F.$

Conversely, let $y^{F} \in [C, D]F//F.$
Then by Lemmas \ref{qinv} and \ref{de},
$$y^{F} \in (c_1^{*}d_1^{*}c_1d_1)^{F} \cdots (c_n^{*}d_n^{*}c_nd_n)^{F} \subseteq ((c_1^{F})^{*}(d_1^{F})^{*})c_1^{F}d_1^{F} \cdots ((c_n^{F})^{*}(d_n^{F})^{*})c_n^{F}d_n^{F} = [c_1^{F}, d_1^{F}] \cdots [c_n^{F}, d_n^{F}]$$
where $c_i \in C$ and $d_i \in D$ for $i = 1, 2, \cdots, n.$
Hence $y^{F} \in [C//F, D//F]$.
It implies that $[C, D]F//F \subseteq [C//F, D//F].$

Therefore $[C//F, D//F] = [C, D]F//F.$
\qed
\medskip

\begin{Corollary}\label{n}

Let $F$ be a normal closed subset of $H$.
Then by Lemma \ref{cq}
$$(H//F)_{s} = \underbrace{[H//F, H//F, \cdots H//F]}_{s} = \underbrace{[H, H, \cdots H]}_{s}F//F = H_{s}F//F,$$ for any non-negative integer s.

\end{Corollary}

\begin{Lemma}\label{com}

$H_{n} = \underbrace{[H, H, \cdots , H]}_{n}$ is a strongly normal closed subset in $H$ for all positive integer $n$.

\end{Lemma}

{\bf Proof.}
It is easy to see that $H_{n} = \underbrace{[H, H, \cdots , H]}_{n}$ is a closed subset of $H$ for all positive integer $n$.
Clearly $H_{1} = H$ is strongly normal in $H$.
Suppose that $H_{i}$ is strongly normal in $H$.
We show that  $H_{i+1} = [H_{i}, H]$ is strongly normal in $H$.
For every $h, y \in H$ and $x \in H_{i}$,
$$h^*[x, y]h = h^*x^*y^*xyh \subseteq h^*x^*hh^*y^*hh^*xhh^*yh = (h^*xh)^*(h^*yh)^*(h^*xh)(h^*yh) = [h^*xh, h^*yh].$$
Since $H_{i}$ is strongly normal in $H$, $h^*xh \in H_{i}$.
And also note that $h^*yh \subseteq H$.
Hence $h^*[x, y]h \in [H_{i}, H] = H_{i+1}$.
It implies from Lemma \ref{de}(ii) that $h^*H_{i+1}h \subseteq H_{i+1}$.
Therefore $H_{i+1}$ is strongly normal in $H$.
As consequence, $H_{n} = \underbrace{[H, H, \cdots , H]}_{n}$ is strongly normal in $H$ for all positive integer $n$.
\qed
\medskip

Let  $S_H(F) =\{ h \in H | h^*Fh \subseteq F \}$.
We call it the strongly normalizer of $F$ in $H$.

\begin{Lemma}\label{sn}

Let $F$ and $K$ are closed subset of $H$ and $K \subseteq F$,
Then $$S_{H//K}(F//K) = S_{H}(F)//K.$$

\end{Lemma}

{\bf Proof.}
Let $h \in H$.
By Lemma 4.3 in \cite{fz1}, Lemma \ref{closed} and Lemma \ref{qinv}, we have
$$h^{K} \in S_{H//K}(F//K)  \Longleftrightarrow$$
$${h^{*}}^{K}F//Kh^{K} \subseteq F//K \Longleftrightarrow$$
$$K^{*}h^{*}FhK \subseteq F \Longleftrightarrow$$
$$h^{*}Fh \subseteq F \Longleftrightarrow$$
$$h \in S_{H}(F) \Longrightarrow$$
$$h^{K} \in S_{H}(F)//K.$$

Conversely, assume that $h^{K} \in S_{H}(F)//K.$
Then $h^{K} = f^{K}$ where $f \in S_{H}(F).$
Hence $h \in KfK \subseteq FfF \subseteq S_{H}(F)$ by Lemma \ref{closed}(i).
It follows from the necessity and sufficiency above that $h^{K} \in S_{H//K}(F//K).$
This completes the proof.

\qed

\textbf{Proof of Theorem \ref{center}}~~

Since $H$ is nilpotent, $H_{n} = \underbrace{[H, H, \cdots , H]}_{n} = [H_{n-1}, H]= 1$ for some non-negative integer $n$.
By Corollary \ref{abel}, $H_{n-1} \subseteq Z(H)$.
By Lemma \ref{com} and Lemma \ref{closed}, we see that $H_{n-1} \leq Z^{*}(H)$.
It implies from Lemma \ref{main1} and Corollary \ref{n} that
$$\underbrace{[H//Z^{*}(H), H//Z^{*}(H), \cdots, H//Z^{*}(H)]}_{n-1} = H_{n-1}Z^{*}(H)//Z^{*}(H) = 1.$$
Hence by Corollary \ref{n} and Lemma \ref{cen},
$$(H//Z^{*}(H))_{n-2} = H_{n-2}Z^{*}(H)//Z^{*}(H) \leq Z^{*}(H//Z^{*}(H)) = Z_{2}^{*}(H)//Z^{*}(H).$$
It implies that $H_{n-2} \leq Z_{2}^{*}(H)$.
As the same as the above proof, we have that $H \leq Z_{n}^{*}(H)$.
Therefore $Z_{n}^{*}(H) = H$.
\qed \\
\medskip

\textbf{Proof of Theorem \ref{ct}}~~

Suppose that $Z_{n}^{*}(H) = H$ for some non-negative integer $n$.
We show that $H//O^{\vartheta}(H)$ is a nilpotent hypergroup.
Firstly, we prove that
$$H_{i} = \underbrace{[H, H, \cdots H]}_{i} \leq Z^{*}_{n-i+1}(H)O^{\vartheta}(H).$$
Since $Z_{n}^{*}(H) = H$, it follows from Lemma \ref{thin residue}(1) that
$$[H//Z_{n-1}^{*}(H), H//Z_{n-1}^{*}(H)] = [Z^{*}(H//Z_{n-1}^{*}(H)), Z^{*}(H//Z_{n-1}^{*}(H))] \leq O^{\vartheta}(H//Z_{n-1}^{*}(H)).$$
Then $[H//Z_{n-1}^{*}(H), H//Z_{n-1}^{*}(H)] \leq O^{\vartheta}(H//Z_{n-1}^{*}(H))=O^{\vartheta}(H)Z_{n-1}^{*}(H)//Z_{n-1}^{*}(H)$ by Lemma \ref{qu}.
Hence $H_2//Z_{n-1}^{*}(H) =[H, H]Z_{n-1}^{*}(H)//Z_{n-1}^{*}(H) = O^{\vartheta}(H)Z_{n-1}^{*}(H)//Z_{n-1}^{*}(H).$
It implies that $H_2 = [H, H] \leq O^{\vartheta}(H)Z_{n-1}^{*}(H).$
Suppose that $H_j=\underbrace{[H, H, \cdots, H]}_{j} \leq O^{\vartheta}(H)Z_{n-j+1}^{*}(H).$
We will show that $H_{j+1}= \underbrace{[H, H, \cdots, H]}_{j+1} \leq O^{\vartheta}(H)Z_{n-j}^{*}(H).$
By Lemma \ref{main1} and Corollary \ref{n},
$$\underbrace{[H//Z_{n-j}^{*}(H), H//Z_{n-j}^{*}(H), \cdots, H//Z_{n-j}^{*}(H)]}_{j+1}$$
$$[\underbrace{[H, H, \cdots, H]}_{j}Z_{n-j}^{*}(H)//Z_{n-j}^{*}(H), H//Z_{n-j}^{*}(H)] $$
$$ \leq [O^{\vartheta}(H)Z_{n-j+1}^{*}(H)//Z_{n-j}^{*}(H), H//Z_{n-j}^{*}(H)]$$
$$= [(O^{\vartheta}(H)Z_{n-j}^{*}(H)//Z_{n-j}^{*}(H))(Z^{*}(H//Z_{n-j}^{*}(H))), H//Z_{n-j}^{*}(H)]$$
$$ \leq O^{\vartheta}(H)Z_{n-j}^{*}(H)//Z_{n-j}^{*}(H)$$
by Lemma \ref{thin residue}.
Now by Lemma \ref{main1} and Corollary \ref{n},
$$ \underbrace{[H//Z_{n-j}^{*}(H), H//Z_{n-j}^{*}(H), \cdots, H//Z_{n-j}^{*}(H)]}_{j+1} = \underbrace{[H, H, \cdots, H]}_{j+1}Z_{n-j}^{*}(H)//Z_{n-j}^{*}(H).$$
It implies that
$$ H_{j+1}=\underbrace{[H, H, \cdots, H]}_{j+1} \leq O^{\vartheta}(H)Z_{n-j}^{*}(H).$$
Note that $O^{\vartheta}(H)Z_{n-j}^{*}(H) = Z_{n-j}^{*}(H)O^{\vartheta}(H)$ by Lemma \ref{closed} and Lemma \ref{main1}.
Hence by induction, we have $H_{i} = \underbrace{[H, H, \cdots H]}_{i} \leq Z^{*}_{n-i+1}(H)O^{\vartheta}(H).$
Now let $i = n+1.$
Then
$$H_{n+1} = \underbrace{[H, H, \cdots, H]}_{n+1} \leq Z^{*}_{0}(H)O^{\vartheta}(H) = O^{\vartheta}(H).$$
It follows that
$$(H//O^{\vartheta}(H))_{n+1} = \underbrace{[H//O^{\vartheta}(H), H//O^{\vartheta}(H), \cdots, H//O^{\vartheta}(H)]}_{n+1}=1. $$
Therefore, $H//O^{\vartheta}(H)$ is nilpotent.
\qed  \\
\medskip

\textbf{Proof of Theorem \ref{ns}}~~
Since $H$ be a finite nilpotent hypergroup,
there exists a closed subset chain
$$ H = H_{1} \supseteq H_{2} = [H, H] \supseteq H_{3} \supseteq \cdots \supseteq H_{n} = 1.$$
By Lemma \ref{com}, we have $H_{i}$ is strong normal in $H_{i-1}$.
Hence by Lemma \ref{cq}, $[H_{i-1}//H_{i}, H//H_{i}] = H_{i}//H_{i} = 1.$
Then by Lemma \ref{cen}, $H_{i-1}//H_{i} \leq Z(H//H_{i}).$
It follows from Lemma \ref{strong} that $H_{i-1}//H_{i}$ is thin.
Then $H_{i-1}//H_{i}$ is an abelian group.
Hence there exits a subgroup chain
$$H_{i}//H_{i} = H_{i1}//H_{i} \leq H_{i2}//H_{i} \leq \cdots \leq H_{im}//H_{i} = H_{i-1}//H_{i}$$
such that $H_{ij}//H_{i}$ is a normal subgroup of $H_{i(j+1)}//H_{i}$ and $(H_{i(j+1)}//H_{i})//(H_{ij}//H_{i})$ is a group of prime order.
It implies that $H_{ij}$ is a strongly normal closed subset of $H_{i(j+1)}$.
And also by \cite[Theorem 4.2]{vz1}, $(H_{i(j+1)}//H_{i})//(H_{ij}//H_{i}) \cong H_{i(j+1)}//H_{ij}$ is a group of prime order.
Therefore $H$ is a finite solvable hypergroup.
\qed
\medskip

\begin{Definition}

A set of hypergroups $\mathfrak{F}$ is said to be {\em a class} if when $H \in \mathfrak{F}$ and $H \cong G$, then $G \in \mathfrak{F}$.
A class is call {\em closed subset hereditary (normal closed subset hereditary, srongly normal closed subset hereditary respectively, (or in short,  $S$-hereditary ($N$-hereditary, $SN$-hereditary))).}
if $F$ is a closed subset $H$ ($F$ is normal in $H$, $F$ is strong normal in $H$) and $H \in \mathfrak{F}$, then $F \in \mathfrak{F}$.
A class is call {\em normal quotient hereditary} (or in short, $NQ$-hereditary) if $F$ is normal in $H$ and $H \in \mathfrak{F}$, it implies that $H//F \in \mathfrak{F}$.

\end{Definition}

\begin{Remark}

 By \cite[Lemma 5.1 and 5.2]{vz1}, the class of solvable hypergroups is $S$-hereditary (of course, $N$-hereditary, $SN$-hereditary) and $NQ$-hereditary.

\end{Remark}

Next we will show that the class of nilpotent hypergroups is also $S$-hereditary (of course, $N$-hereditary, $SN$-hereditary) and $NQ$-hereditary.

\begin{Proposition}\label{s}
Let $H$ be a nilpotent hypergroup and $F$ is a closed subet of $H$.
Then $F$ is nilpotent.
\end{Proposition}

{\bf Proof.}
Since $H$ be a nilpotent hypergroup,
$H_n = \underbrace{[H, H, \cdots, H]}_{n} = 1$ for some positive integer.
It follows from $F$ is a closed subet of $H$ that
 $$F_{n} = \underbrace{[F, F, \cdots, F]}_{n} \leq H_n = \underbrace{[H, H, \cdots, H]}_{n} = 1.$$
Therefore $F$ is nilpotent.
\qed
\medskip

\begin{Proposition}\label{NQ}
Let $H$ be a nilpotent hypergroup and $F$ is a normal closed subet of $H$.
Then $H//F$ is nilpotent.
\end{Proposition}

{\bf Proof.}
Since $H$ be a nilpotent hypergroup,
$H_n = [H, H, \cdots, H] = 1$ for some positive integer.
Hence by Corollary \ref{n},
$$(H//F)_n=\underbrace{[H//F, H//F, \cdots, H//F]}_{n} = H_nF//F ={1}.$$
Then $H//F$ is nilpotent.
\qed

\medskip

We call a closed subset $M$ is {\em maximal} in $H$ if $M \neq H$ and when there exits a closed subset $K$ of $H$ such that $M \subseteq K \subset H$, then $M = K$.

\vspace{0.25cm}

\textbf{Proof of Theorem \ref{strongly}}~~

We will divide the proof into the following steps.

{\bf Step $(1)$} {\sl If $F$ is a proper closed subset of $H$, then $F$ is a proper closed subset of $S_H(F)$.}

Let $F$ is a proper closed subset of $H$.
Since $H$ is nilpotent, there exits some $i$, such that $F \supseteq H_{i+1}$ and $F \nsupseteq H_{i}$.
By Lemmas \ref{cq} and \ref{com}, we have $[H_{i}//H_{i+1}, F//H_{i+1}] = [H_{i}, F]H_{i+1}//H_{i+1} =1$.
Thus by Lemma \ref{cen},
$$H_{i}//H_{i+1} \subseteq C_{H//H_{i+1}}(F//H_{i+1}) \subseteq S_{H//H_{i+1}}(F//H_{i+1}).$$
It follows from Lemma \ref{sn} that
$$H_{i}//H_{i+1} \subseteq S_{H//H_{i+1}}(F//H_{i+1}) = S_{H}(F)//H_{i+1}.$$
Hence $H_{i} \subseteq S_{H}(F)$.
It implies that $F$ is a proper closed subset of $S_H(F)$.

{\bf Step $(2)$} {\sl If $M$ is a maximal closed subset of $H$, then $M$ is strongly normal in $H$.}

It follows from Step $(1)$ that $M$ is a proper closed subset of $S_H(M)$.
Note that $S_H(M)$ is a closed subset of $H$.
In fact, for every $a, b \in S_H(M)$,
$${(a^{*}b)}^{*}M(a^{*}b) = b^{*}aMa^{*}b \subseteq  b^{*}Mb \subseteq M.$$
So $a^{*}b \in S_H(M)$.
Hence $S_H(M)$ is a closed subset of $H$.
But since $M$ is a maximal closed subset of $H$, $S_H(M) = H$.
Hence $M$ is strongly normal in $H$.

{\bf Step $(3)$} {\sl For every closed subset $E$ of $H$ and $E \neq 1$, $E$ is strongly subnormal in $H$.}

In the case that $E=H$, clearly $E$ is strongly subnormal in $H$.
Now assume that $E \neq H$ is a proper closed subset of $H$.
Then there exits a closed subset chain
$$E = E_{0} \leq E_{1} \leq \cdots \leq E_{n} = H$$
such that $E_{i}$ is a maximal closed subset of $E_{i+1}$ for $i=0, 1, \cdots, n-1$.
By Proposition \ref{s}, $E_{i}$ is nilpotent for $i=0, 1, \cdots, n$.
It follow from  Step $(2)$ that $E_{i}$ is strongly normal in $E_{i+1}$ for $i=0, 1, \cdots, n-1$.
Therefore $E$ is strongly subnormal in $H$.
This completes the proof.
\qed

\section{Some remarks and  questions}

In the theory of finite group, nilpotent groups play a very important role.
There are numerous research results on nilpotent groups (see monographs \cite{GuoI,Rob,H}).
The concept of nilpotent hypergroups mentioned in this paper is actually a generalization of classical nilpotent groups.
In fact, nilpotent groups can be regarded as tin nilpotent hypergroups and also tin nilpotent hypergroups can be regarded as nilpotent groups.
We not only generalized the concept of nilpotent groups, but also extended many classic results of nilpotent groups.
Below, we will provide some reminders and some applications of my obtained results.

\begin{Remark}\label{z}

When $H$ is a group,
then $Z^{*}(H) = Z(H)$ is a normal subgroup of $H$ and $O^{\vartheta}(H) = 1$.
And  $Z_{\infty}^*(H) = Z_{\infty}(H)$ is a hypercenter of $H$.
Moreover, the chain in Definition \ref{hyper} is a upper centre subgroup chain in a nilpotent group.

\end{Remark}

Combining Remark \ref{z}, we can derive the following result from Theorems \ref{center} and \ref{ct}.

\begin{Corollary}

A group $G$ is nilpotent if and only if $Z_{n}(G)=G$ for some positive integer $n$.

\end{Corollary}

From Theorem \ref{strongly}, we can get the following Corollary.

\begin{Corollary}

Let $H$ be a finite nilpotent group.
Then every subgroup of $H$ is subnormal in $H$.

\end{Corollary}

\begin{Remark}

As we known, if every subgroup of a finite group $H$ is subnormal in $H$, then $H$ is nilpotent.
However, we do not know if this conclusion holds for hypergroups.

\end{Remark}

Then we have the following question.

\begin{Question}

Let $H$ be a finite nilpotent hypergroup.
Is any closed subset $F$ of $H$ is subnormal in $H$ $?$

\end{Question}

Following \cite{b1}, a hypergroup $H$ is called residually thin (or in short $RT$)
if there exists a chain of closed subsets
$\{1\} = F_0 \subset F_1 \subset \ldots \subset F_n = H$ such that $F_{i}//F_{i-1}$ is thin for all $1 \leq i \leq n$
and the valency of a finite $RT$ hypergroup $H$ is the integer
$$n_{H} = \prod_{i=1}^{n} |F_{i}//F_{i-1}|.$$
Any closed subset $C$ with that $n_C$ is a $p$-number is called a $p$-subset;
A Sylow $p$-subset of $H$ is a $p$-subset $C$ such that $n_{H}/n_{C}$ is a $p'$-number.

We know that every Sylow subgroup of a finite nilpotent group $G$ is normal in $G$.
However we do not know if  this conclusion holds for finite $RT$ nilpotent hypergroups.
Hence we have the following question.

\begin{Question}

Let $H$ be a finite $RT$ nilpotent hypergroup.
Is any Sylow closed subset $P$ of $H$ is normal in $H$ $?$

\end{Question}


\begin{thebibliography}{99}

\bibitem{b1}
H. Blau, P.-H. Zieschang, Sylow theory for table algebras, fusion rule algebras, and hypergroups, J.
Algebra \textbf{273} (2004) 551-570.

\bibitem{fz1}
C. French, P.-H. Zieschang, On residually thin hypergroups, J. Algebra \textbf{551} (2020) 93-118.

\bibitem{GuoI}
W. Guo, \textit{The Theory of Classes of Groups}, Science Press, Kluwer Academic Publishers, Beijing-New York-Dordrecht-Boston-London, 2000.

\bibitem{H}
B. Huppert, \emph{Endliche Gruppen I}, Springer-Verlag, Berlin, 1967.

\bibitem{m1}
F. Marty, Sur une generalization de la notion de groupe, in: 8th Congress Math. Scandinaves, Stockholm, 1934, pp. 45-49.

\bibitem{Rob}
D. J. S. Robinson, A Course in the Theory of Groups, Springer-Verlag, New York, 1982.

\bibitem{tz1}
R. Tanaka, P.-H. Zieschang, On a class of wreath products of hypergroups and association schemes,
J. Algebraic Comb. \textbf{37} (2013) 601-619.

\bibitem{vz1}
A. Vasil'ev, P.-H. Zieschang, Solvable hypergroups and a generalization of Hall's theorems on solvable groups to association schemes, J. Algebra \textbf{594} (2022) 733-750.

\bibitem{z2}
P.-H. Zieschang, Hypergroups all non-identity elements of which are involutions, in: Advances in
Algebra 305-322, in: Springer Proc. Math. Stat., vol. 277, Springer, Cham, 2019.

\end{thebibliography}
\end{document}